\documentclass[11pt]{article}
\usepackage{amsmath,amsfonts,amssymb,amscd,amsthm}

\usepackage{graphicx}

\newtheorem{theo}{Theorem}
\newtheorem{lem}{Lemma}[section]
\newtheorem{prop}{Proposition}[section]

\newtheorem{cor}{Corollary}[section]
\newtheorem{rem}{Remark}[section]

\makeatletter \@addtoreset{equation}{section} \makeatother

\newcommand{\mC}{\mathbb{C}}

\newcommand{\mR}{\mathbb{R}}

\newcommand{\bA}{{\bf A}}

\newcommand{\bQ}{{\bf Q}}
\newcommand{\bS}{{\bf S}}

\newcommand{\bd}{{\bf d}}

\newcommand{\bn}{{\bf n}}

\newcommand{\bL}{\mathbf{L}}

\newcommand{\bv}{{\bf v}}

\newcommand{\calD}{{\cal D}}

\newcommand{\calF}{{\cal F}}

\newcommand{\calL}{{\cal L}}

\newcommand{\calO}{{\cal O}}
\newcommand{\calP}{{\cal P}}

\newcommand{\eps}{\varepsilon}

\newcommand{\diag}{\operatorname{diag}}

\newcommand{\im}{\operatorname{Im}}

\newcommand{\re}{\operatorname{Re}}

\newcommand{\spec}{\operatorname{Spec}}

\begin{document}

\date{}

\title
{On the problem of stability of viscous shocks}
\author{ S.~Bolotin and D.~Treschev\\
Steklov Mathematical Institute of Russian Academy of Sciences}

\maketitle

\centerline{Dedicated to the memory of Alexey Borisov}
 
\begin{abstract}

We consider the problem of spectral stability of traveling wave solutions $u=\gamma(x-Wt)$
for a system of viscous conservation laws $\partial_t u + \partial_x F(u) = \partial^2_x u$.
Such solutions correspond to heteroclinic trajectories $\gamma$ of a system of ODE.
In general conditions of stability can be obtained only numerically.
We propose a model class of piece-wise linear (discontinuous) vector fields $F$ for which the stability
problem is reduced to a linear algebra problem. We show that the stability problem makes sense in such low regularity and construct
several examples of stability loss. Every such example can be smoothed to provide a smooth example of the same phenomenon.

\end{abstract}

\section{Introduction}

Consider the following system of nonlinear parabolic PDEs
\begin{equation}
\label{pde}
  \partial_t u + \partial_x F(u) = \nu\partial^2_x u,\qquad
  t,x\in\mR, \quad u=u(x,t)\in\mR^n ,\quad \nu>0,
\end{equation}
where $F:\mR^n\to\mR^n$ is a smooth vector field. System (\ref{pde}) is called a system of viscous conservation laws.
It is a viscous regularization of the system
\begin{equation}
\label{eq:hyp}
\partial_t u + \partial_x F(u) = 0.
\end{equation}
This system is said to be strictly hyperbolic if for all $u$ the Jacobi matrix of $F$
has $n$ real distinct eigenvalues. Solutions of system (\ref{eq:hyp}) develop singularities, usually called shocks.
For small viscosity $\nu$, solutions of system (\ref{pde})  are expected to represent the structure of the shock
in a narrow band  around the singularity which tends to zero as $\nu\to 0$.
However, scaling the variables, one can make $\nu=1$. We will assume this from now on.

A solution of (\ref{pde}) of the form
\begin{equation}
\label{tws}
  u(x,t) = \gamma(\xi), \quad
  \xi = x - Wt, \quad
  \gamma(\pm\infty)=\lim_{\xi\to\pm\infty} \gamma(\xi) = u^\pm,
\end{equation}
is called a traveling wave solution (TWS), or also a viscous shock.
Then there exists a constant $D\in\mR^n$ such that $\gamma:\mR\to\mR^n$
is a heteroclinic solution of the ODE system
\begin{equation}
\label{ode}
  u' = f(u), \qquad
  f(u) = F(u) - Wu + D, \quad
  u' = du/d\xi,
\end{equation}
which goes from $u^-$ to $u^+$.
Hence, $u^\pm$ are equilibrium (critical) points  of (\ref{ode}): $f(u^\pm)= 0$.

Let $u=\gamma(\xi)$ be such a heteroclinic.
We are interested in the problem of stability of the corresponding TWS  of system (\ref{pde})
(stability means orbital stability:  modulo the phase shift $\xi\mapsto\xi+c$).
If we make a change of variables $x\mapsto \xi=x-W t$ (pass to a moving frame), system (\ref{pde}) becomes
\begin{equation}
\label{eq:pde}
  \partial_t u=\calF(u),\qquad \calF(u)= - \partial_\xi f(u)+\partial^2_\xi u,
\end{equation}
and $u=\gamma(\xi)$ is an equilibrium (stationary) solution: $\calF(\gamma)=0$.

In this paper we deal only with stability for the linearized system which is also called the spectral stability.
Substituting
$$
u = \gamma(\xi) +  v(\xi,t)
$$
into (\ref{eq:pde}) ($v$ is assumed to be small), in the first approximation in $v$ we obtain the linearized system
$$
\partial_t v=\calL v,
$$
where $\calL$ is the linear operator (linearization of $\calF$):
\begin{equation}
\label{eq:calL}
\calL v=\big(v' -  \partial f(\gamma(\xi)) v\big)'=\partial_\xi\big( \partial_\xi v-\partial f(\gamma(\xi))v \big).
\end{equation}

Spectral stability of the TWS (\ref{tws}) means that there is no spectrum of $\calL$ in the right half-plane:
 $\spec\calL\subset \{\lambda\in\mC : \re\lambda \le 0\}$.

\begin{rem}
Usually $\calL$ is regarded as an operator on $L^2$. In general it has essential spectrum.
However, under mild assumptions, there is no essential spectrum in the right half plane. Moreover, sometimes it is possible to prove that there is a gap between the essential spectrum and the imaginary axis, see \cite{Gar-Zum}. Hence for stability problems only eigenvalues of $\calL$ matter.
\end{rem}

\begin{rem}
\label{rem:lam=0}
$v=\gamma'$ is a solution of the variational equation $v'=\partial f(\gamma)v$. Hence
$\lambda=0$ is always an eigenvalue, and the corresponding eigenfunction is $v=\gamma'$.
\end{rem}

It is convenient to rewrite the eigenvalue problem $\calL v=\lambda v$ as a first order spectral problem
\begin{equation}
\label{eq:system}
v'=\partial f(\gamma)v+z, \quad z'=\lambda v.
\end{equation}
We shortly write this system as
\begin{equation}
\label{eq:bL}
\bL(\lambda)\bv = \bv'-\bA(\xi,\lambda)\bv=0,\qquad
  \bv(\pm\infty)=0,
\end{equation}
where
$$
\bv = \begin{pmatrix} v\\ z \end{pmatrix} \quad \mbox{and}\quad
    \bA(\xi,\lambda)
  = \begin{pmatrix} \partial f(\gamma(\xi)) & I \\ \lambda I & 0 \end{pmatrix}.
$$
Let
$$
  E^\pm(\lambda)=\{\bv(0): \bL(\lambda)\bv=0,\;\bv(\pm\infty)=0\}.
$$
Then $\lambda$, $\re\lambda\ge 0$, is an eigenvalue of the spectral problem (\ref{eq:bL}) iff
\begin{equation}
\label{eq:capE1}
E^+(\lambda)\cap E^-(\lambda)\ne \{0\}.
\end{equation}

For $n=1$ the problem of stability of TWS's is well studied.
The first results were based on the maximum principle \cite{IO}.
They deal with generalizations of the Burgers equation and establish nonlinear stability.
Other methods were used in \cite{MN,Satt}.

For $n\ge 2$ the main tool to study spectral stability is the Evans function  \cite{Evans}.
The Evans function is an analog of the characteristic polynomial for infinite dimensions:
eigenvalues of $\calL$ in the right half plane $\re\lambda\ge 0$ are zeros of the Evans function.
In many papers, spectral stability was established numerically with the help of a
combination of the Evans function method and the argument principle, see e.g.  \cite{Chug-Pol,Ilich}.
Spectral instability was rigorously established  in several papers, see e.g.\  \cite{Chug-Tre,Gar-Zum},
but only for very special vector fields.
General results on passage from linear (spectral)  to nonlinear stability were obtained
in e.g.\ \cite{Good,LZ1,LZ2,Liu,How-Zum,Zum-How},
see also a comprehensive survey in \cite{Zum:Survey}.
Under mild natural assumptions it was proved that spectral stability of a shock
(no roots of the Evans function in the right half plane $\re\lambda\ge 0$
except zero eigenvalue with minimal order  -- the dimension of the family of shocks)
implies nonlinear stability --  convergence of perturbed solutions to a viscous shock in the same family.
Spectral instability (roots of the Evans function with $\re\lambda>0$) implies nonlinear instability,
exactly as in the classical Lyapunov theorem.

Most rigorous instability results for concrete systems were obtained for $n=2$ and
quadratic vector fields, see e.g.\ \cite{Gar-Zum,Chug-Tre} and the survey  \cite{Zum:Survey}.
Then there are simplifications: for example heteroclinics connecting saddles
(they correspond to undercompressive shocks) are straight segments \cite{Chicone}.

However, it seems that mechanisms of spectral instability for TWS's still remain unclear.
In this paper we propose a model class of vector fields (piece-wise linear fields)
for which explicit stability results may be obtained via elementary linear algebra.
We hope that understanding the phenomena
of stability loss for this model class will give us more insight in what happens in general.

We start with the situation when the vector field $f$ is linear on two half-spaces separated by a hyperplane $\Sigma\subset\mR^n$.
Then $f$ is discontinuous on $\Sigma$. We reduce the problem of stability of viscous shocks,
corresponding to heteroclinic solutions of (\ref{ode}), to a  linear algebra problem.
For $n=2$ we present simple sufficient conditions for spectral stability and instability.

Then we give an example of stability loss  for overcompressive shocks corresponding to
a heteroclinic connecting an unstable node $u^-$ with a stable node $u^+$.
In this case there is a one-parameter family of heteroclinics $\gamma_s$ joining $u^-$ and $u^+$,
see Figure \ref{fig:Over}.
We give examples where:

\begin{figure}
\label{fig:Over}
\includegraphics[scale=0.9]{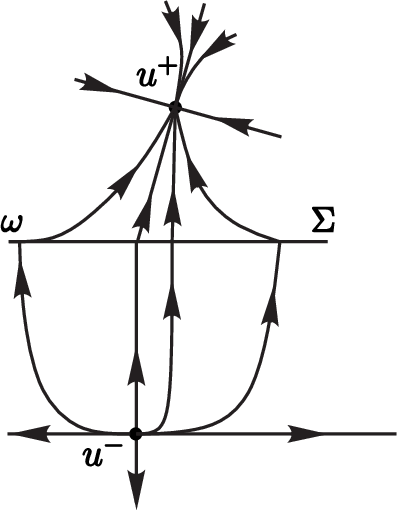}
\caption{Overcompressive shock}
\end{figure}

\begin{itemize}
\item  TWSs corresponding to all heteroclinics $\gamma_s$ are spectrally stable;
\item  TWS corresponding to all heteroclinics $\gamma_s$ are spectrally unstable;
\item  There is stability loss  inside this family:
a pair of complex conjugate eigenvalues of the spectral problem crosses the imaginary axis.
More precisely, there exist $s_0<s_1<s_2$ such that for $s_0<s<s_1$ the TWS corresponding to $\gamma_s$
is spectrally stable, and for $s=s_1$ two conjugate eigenvalues $\lambda(s)\ne \bar\lambda(s)$
cross the imaginary axis: $\re\lambda(s)<0$ for $s_1-\delta<s<s_1$ and $\re\lambda(s_1)=0$, $\im\lambda(s_1)\ne 0$.
Thus the TWS looses spectral stability at $s=s_1$
through an Andronov-Hopf type bifurcation.
For large $s>s_2$ a real positive eigenvalue $\lambda\sim s^{-1}$ exists.
\end{itemize}

The most technical example we consider is when  instability of  a Lax shock is associated with the bifurcation
of a passage through a saddle connection.
Such a bifurcation appears in many models.
One of the simplest examples may be found in \cite{Chug-Tre}.

Suppose that the vector field $f$ depends on the parameter $\eps$:
$$
  f = f^\eps, \qquad
  -\eps_0 < \eps < \eps_0.
$$
We assume that for $\eps=0$ system (\ref{ode}) has 3 equilibrium points $u^{-}$, $u^{\times}$, $u^{+}$
such that

\begin{figure}
\includegraphics[scale=0.8]{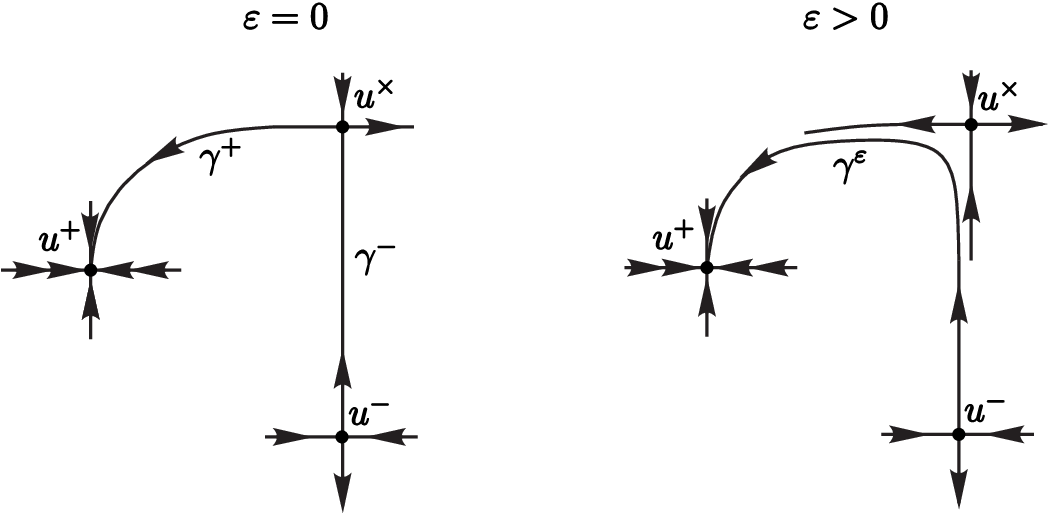}
\caption{Splitting of a saddle connection}
\label{fig:Bifur}
\end{figure}

\begin{itemize}
\item
$u^{-}$ and $u^{\times}$ are saddles while $u^{+}$ is a stable node,

\item
there are heteroclinics $\gamma^-$ from $u^{-}$ to $u^{\times}$ and
$\gamma^+$ from $u^{\times}$ to $u^{+}$,
see Figure \ref{fig:Bifur} left.
\end{itemize}

Then for small $|\eps|$ the vector field $f^\eps$ has equilibrium points
$u^{-\eps}$, $u^{\times\eps}$, $u^{+\eps}$ smoothly depending on $\eps$.
The saddle connection $\gamma^-$ generically breaks and a heteroclinic solution
$\gamma^\eps$ from $u^{-\eps}$ to $u^{+\eps}$ may appear
for $\eps>0$ or $\eps<0$, see Figure \ref{fig:Bifur} right. This is determined by an appropriate Melnikov function.  We will show that the following situation may occur:

\begin{itemize}

\item for $\eps=0$ TWS's corresponding to $\gamma^\pm$ are linearly stable;

\item for any small $\eps>0$ the eigenvalue problem associated with the TWS $\gamma^\eps$
has an  eigenvalue $\lambda = c\eps + o(\eps)$, where $c>0$ is a constant.
Hence the TWS $\gamma^\eps$ is linearly unstable.
\end{itemize}

For $n=2$ a heteroclinic from a saddle to a stable node corresponds to a Lax shock
for the hyperbolic system (\ref{eq:hyp}),
while a heteroclinic from a saddle to a saddle corresponds  to an undercompressive shock.
So, physically speaking we have a bifurcation of a stable Lax shock to an unstable one
through merging with a stable undercompressive shock.

We work in the class of piece-wise linear vector fields, having breaks along singularity lines.
So our first goal will be to show that the spectral stability problem makes sense in such low regularity.
Once we establish spectral instability  of TWS's for this class of vector fields,
smoothing will provide smooth examples of the same phenomena.

In Section \ref{sec:reg} we discuss regularization of piece-wise continuous vector fields
and, as a consequence, derive the spectral problem for heteroclinics of such vector fields.
The proof is postponed to Section \ref{sec:App}.
In Section \ref{sec:examples} we give precise description of the examples presented above
and formulate three main theorems. In Section \ref{sec:simple} we discuss general properties of
the spectral problem for heteroclinics of piece-wise linear vector fields with a single discontinuity
and then prove Theorems \ref{theo:simple} and \ref{theo:over}.
Theorem \ref{theo:main} is proved in Section \ref{sec:proof}.

\section{Regularization of discontinuities}

\label{sec:reg}

We need to make sense of the operator (\ref{eq:calL}) and the corresponding spectral problem (\ref{eq:system})
for discontinuous vector fields $f:\mR^n\to\mR^n$.
Suppose the vector field $f$ has a break on the
hyperplane\footnote{The results hold with almost no difference when $\Sigma$ is a smooth hypersurface.
However, we will not need this generalization.}  $\Sigma\subset \mR^n$ dividing $\mR^n$ in
two open regions $U^-$ and $U^+$. More precisely, suppose
$$
    f(u)
  = \left\{\begin{array}{cc} f^+(u), &  u\in U^+,  \\
                             f^-(u), &  u\in U^-,
           \end{array}\right.
$$
where $f^\pm$ are smooth in $\overline{U^\pm}=U^\pm\cup\Sigma$. Let $\bn$ be the unit normal to $\Sigma$
in the direction from $U^-$ to $U^+$.
For trajectories of $u'=f(u)$ crossing  $\Sigma$ from $U^-$ to $U^+$ at a point $w\in\Sigma$,  we have
\begin{equation}
\label{eq:trans}
  \langle f^\pm(w),\bn\rangle > 0.
\end{equation}
Inequality (\ref{eq:trans}) implies the existence and uniqueness of solutions:
they cross $\Sigma$ transversely.

Let $\gamma:[-T,T]\to \mR^n$ be a trajectory crossing $\Sigma$ at a point $\gamma(0)=w$.
We need to make sense of the eigenvalue equation (\ref{eq:calL}), or the spectral system (\ref{eq:system}),
when $f(\gamma(\xi))$ has a break at $\xi=0$.
A natural way is to change $f$ in a narrow layer near the singularity replacing it by a  smooth vector field
$f^\mu$, and then pass to the limit $\mu\to 0$.

We may assume that $\Sigma$ passes through 0:
$$
\Sigma=\{u:\langle \bn,u\rangle=0\},\quad   U^\pm=\{u:\pm \langle \bn,u\rangle>0\}.
$$

 Let $\phi\in C^\infty(\mR)$ be a function such that
$$
    \phi(y)
  = \left\{\begin{array}{cc} 1, & y\ge 1  \\
                             0, & y\le 0
           \end{array}
    \right., \qquad
    \phi'(y) > 0, \quad
    0<y<1.
$$
Take small $\mu>0$ and  define the vector field $f^\mu$ by
\begin{equation}
\label{eq:phi}
    f^\mu(u)
  = \phi(\langle \bn,u\rangle/\mu) f^+(u) + (1-\phi(\langle \bn,u\rangle/\mu)) f^-(u).
\end{equation}
Then $f^\mu$ is smooth and $f^\mu=f$ except in the layer $U_\mu=\{u:0<\langle \bn,u\rangle<\mu\}$.
As $\mu\to 0$, $f^\mu(u)\to f(u)$ and $\partial f^\mu(u)\to \partial f(u)$ pointwise for $u\notin\Sigma$.

Let $\gamma^\mu:[-T,T]\to\mR^n$ be a solution of the regularized ODE
$$
 \gamma^{\mu\,\prime} = f^\mu(\gamma^\mu), \qquad
  \gamma^\mu(0) = w.
$$
By the transversality assumption (\ref{eq:trans}), as $\mu\to 0$,
$\gamma^\mu(\xi)$ converges pointwise on $[-T,T]$ to the piece-wise smooth solution $\gamma(\xi)$ of $u'=f(u)$
with $\gamma(0)=w$.

We regularize the eigenvalue  problem (\ref{eq:system}) for the trajectory $\gamma$, replacing $\partial f(\gamma(\xi))$
by $\partial f^\mu(\gamma^\mu(\xi))$:
\begin{equation}
\label{eq:zv}
 v'=\partial  f^\mu(\gamma^\mu(\xi))v+z,\quad z'=\lambda v.
\end{equation}

Denote
$$
\hat f^\pm = f^\pm(w),\quad \Delta f=\hat f^+-\hat f^-.
$$
Let $v^\mu(\xi)$, $z^\mu(\xi)$ be a solution of (\ref{eq:zv}) with the initial condition $v^\mu(0)=v^-$, $z^\mu(0)=z^0$.
The next proposition is proved in Section \ref{sec:smoothing}.

\begin{prop}\label{prop:smooth}
As $\mu\to 0$, $z^\mu(\xi)$ converges uniformly on $[-T,T]$ to $z(\xi)$, while $v^\mu(\xi)$ is uniformly bounded
and   converges pointwise on  $[-T,T]\setminus \{0\}$ to  $v(\xi)$.
The functions $z(\xi),v(\xi)$ are smooth on $[-T,0]$ and $[0,T]$ and are solutions of the initial value problems
$$
\begin{array}{lllllll}
 v'=\partial f^-(\gamma(\xi))v+z, \quad & z'=\lambda v,\quad    &-T\le \xi\le 0 ,
 \quad &v(-0)=v^-,\; & z(-0)=z^0,\\
v'=\partial  f^+(\gamma(\xi))v+z, \quad &   z'=\lambda v,\quad    &0\le \xi\le T ,
  \quad &v(+0)=v^+, \; &z(+0)=z^0,
  \end{array}
$$
where
\begin{equation}
\label{eq:S}
v^+=Sv^-:=v^-+\Delta f \frac{\langle v^-,\bn\rangle}{\langle \hat f^-,\bn\rangle}.
\end{equation}
\end{prop}

\begin{rem}
The correspondence between $v^-$ and $v^+$ is symmetric:
$$
v^-=v^+-\Delta f \frac{\langle v^+,\bn\rangle}{\langle \hat f^+,\bn\rangle}.
$$
\end{rem}

The limit function $v(\xi)$  is a smooth solution of $\calL^\pm v=(v'-\partial f^\pm(\gamma(\xi))v)'=\lambda v$
on  $[0,T]$  and $[-T,0]$ respectively.
It is discontinuous at $\xi=0$:
\begin{eqnarray*}
&  v^- = \lim_{\xi\nearrow 0} \lim_{\mu\searrow 0} v^\mu(\xi)\ne
   v^+ = \lim_{\xi\searrow 0} \lim_{\mu\searrow 0} v^\mu(\xi), & \\
&  {v^-}' = \lim_{\xi\nearrow 0} \lim_{\mu\searrow 0} {v^\mu}'(\xi)\ne
   {v^+}' = \lim_{\xi\searrow 0} \lim_{\mu\searrow 0} {v^\mu}'(\xi). &
\end{eqnarray*}

We will need  the formula for the operator $S$  in (\ref{eq:S}) for $n=2$ and $\Sigma\subset\mR^2$ a line.

\begin{cor}
\label{cor:smoothing}
If $\Sigma=\{u_2=0\}$  and $\bn=(0,1)$, then
\begin{equation}
\label{eq:S-}
     S=
     \begin{pmatrix} 1  &  (\hat f_1^+ - \hat f_1^-) / \hat f_2^- \\
                             0  &   \hat f_2^+ / \hat f_2^-
           \end{pmatrix}.
\end{equation}
If  $\Sigma=\{u_1=0\}$ and $\bn=(-1,0)$, then
\begin{equation}
\label{eq:S+}
     S=
     \begin{pmatrix} \hat f_1^+ / \hat f_1^-  & 0   \\
                             (\hat f_2^+ - \hat f_2^-) / \hat f_1^-  & 1
           \end{pmatrix}.
\end{equation}

\end{cor}

Suppose now that $\gamma:\mR\to\mR^n$ is a heteroclinic joining $u^\pm\in U^\pm$ and intersecting the break surface $\Sigma$ once at $\xi=0$.
We obtain that the correct way to define the spectral system $\bL(\lambda)\bv=0$ in (\ref{eq:system})
for a piece-wise continuous vector field $f$ is
\begin{equation}
\label{eq:L+-}
     \left\{\begin{array}{lll}
        z'=\lambda v,\quad v'=\partial f^-(\gamma(\xi))v+z,   &\xi <0 , \quad &v(\pm\infty)=z(\pm\infty)=0, \\

        z'=\lambda v,\quad v'=\partial f^+(\gamma(\xi))v+z,   &\xi >0 , \quad &v(+0)=Sv(-0),\;  z(+0)=z(-0).
                   \end{array}\right.
\end{equation}
Let $\bL^\pm(\lambda)$ be the operators (\ref{eq:bL}) defined by the vector fields $f^\pm$. Set
\begin{equation}
\label{eq:Epm}
E^\pm(\lambda)=\{\bv(0):\bL^\pm(\lambda)\bv= 0,\; \pm\xi>0,\; \bv(\pm\infty)=0\}.
\end{equation}
Then $\lambda$ is an eigenvalue of the spectral problem
associated to $\gamma$ iff
\begin{equation}
\label{eq:capE}
E^+(\lambda)\cap \bS E^-(\lambda)\ne\{0\},\qquad
\bS=\begin{pmatrix}
S & 0\\
0 & I
\end{pmatrix}.
\end{equation}
If there is no break of $f$, then $S=I$ and we obtain (\ref{eq:capE1}).

Suppose the equilibria  $u^\pm$ of the systems $u'=f^\pm(u)$ are hyperbolic.
Let $W^{s,u}(u^\pm)\subset \overline{U^\pm}$ be the (local) stable and unstable manifolds
and $W_\mu^{s,u}(u^\pm)$ the (local) stable and unstable manifolds for the smoothed system $u'=f^\mu(u)$.
Then we have:

\begin{prop}
$$
  T_{w}W_\mu^u(u^-)= T_{w} W^u(u^-), \quad
  \lim_{\mu\to 0} T_{w}W_\mu^s(u^+)= S^{-1}T_{w}W^s(u^+),\qquad w = \gamma(0)\in \Sigma.
$$
\end{prop}

The first formula is evident: $f^\mu=f^-$ in $U^-\cup\Sigma$.
To prove the second formula, notice that for $\lambda=0$ the eigenvalue problem (\ref{eq:zv})
degenerates into the variational equation $v'=\partial f^\mu(\gamma) v$.
Hence the proof of Proposition \ref{prop:smooth} in Section \ref{sec:App} works.
\qed

\begin{cor}
If the transversality condition holds:
\begin{equation}
\label{eq:transv}
  S T_{w}W^u(u^-) + T_{w}W^s(u^+)=\mR^n,
\end{equation}
then for small $\mu>0$ there exists  a smooth family  $\gamma^\mu$ of heteroclinics of
$u'=f^\mu(u)$ connecting $u^-$ and $u^+$ such that $\gamma^\mu\to \gamma$ pointwise as $\mu\to 0$.
\end{cor}

Indeed,  transverse intersections survive small perturbations of the manifolds.
\qed

\medskip

If moreover $\gamma$ corresponds to a Lax shock, i.e.\ the compressivity index
$$
d=\dim W^u(u^-)+\dim W^s(u^+)- n=1,
$$
then
$$
ST_{w}W^u(u^-)\cap T_{w}W^s(u^+)=\mR \hat f^+,
$$
and there is a unique (modulo the shift $\xi\mapsto\xi+c$) one-parameter family $\gamma^\mu$ of heteroclinics.
For overcompressive shocks with $d>1$,
there are many such families.

If $\gamma$ corresponds to an undercompressive shock, i.e.\ $d\le 0$,
then (\ref{eq:transv}) never holds.
Then in general such a family does not exist:
this nongeneric case is not stable under perturbations.
However, undercompressive shocks usually appear when the vector field $f$ has some symmetry,
for example invariant under a reflection, and $\gamma$ is preserved by the reflection.
Then, if smoothing of $f$ respects the symmetry,
the vector field $f^\mu$ still possesses  a heteroclinic.

Assuming that the family $\gamma^\mu$ of heteroclinics exists,
let $\bL^\mu(\lambda)=0$ be the corresponding eigenvalue problem (\ref{eq:bL}).

\begin{prop}
Let $\lambda$, $\re\lambda>0$, be a simple eigenvalue of the eigenvalue problem $\bL(\lambda)\bv=0$
defined by $\gamma$.
For small $\mu>0$ the eigenvalue problem $\bL^\mu(\lambda)\bv=0$
has a simple eigenvalue $\lambda^\mu=\lambda+O(\mu)$.
\end{prop}

Indeed, $\lambda$  is a simple eigenvalue if the spaces
$\bS E^-(\lambda)$ and $E^+(\lambda)$
have one-dimensional intersection and the determinant $\Delta(\lambda)$ of the corresponding system
of $2n$ linear equations has a simple zero. But under conditions above, the perturbed spaces
$E_\mu^\pm(\lambda)$  satisfy
$$
E_\mu^-(\lambda)=E^-(\lambda),\quad \lim_{\mu\to 0}E_\mu^+(\lambda)=\bS^{-1}E^+(\lambda).
$$
Hence for small $\mu>0$, the corresponding determinant $\Delta_\mu(\lambda)=\Delta(\lambda)+O(\mu)$
has a simple zero $\lambda_\mu=\lambda+O(\mu)$.
\qed

\section{Examples}

\label{sec:examples}

All examples we give are for piece-wise linear discontinuous vector fields $f:\mR^2\to\mR^2$.
As explained in the previous section, every such field can be smoothed
to provide a smooth example of the same phenomenon.

First consider simplest possible shocks. Let
\begin{equation}
\label{u+-}
  u^- = \begin{pmatrix} u_1^- \\ u_2^- \end{pmatrix}, \quad
  u^+ = \begin{pmatrix} u_1^+ \\ u_2^+ \end{pmatrix}, \qquad
  \pm u_2^\pm > 0.
\end{equation}
Set $\Sigma = \{u_2=0\}$ and
\begin{equation}
\label{eq:f(u)}
\begin{array}{lll}
&    f(u)
  =  \left\{ \begin{array}{lcl}
               Q^-(u-u^-)\quad                 & \mbox{ if } & \mbox{ $u_2 < 0$}, \\
               Q^+(u-u^+)\quad                 & \mbox{ if } & \mbox{ $u_2 > 0$}.
            \end{array}
     \right. &
\end{array}
\end{equation}
The matrices   $Q^\pm = \diag(h_1^\pm,h_2^\pm)$ are diagonal:, where at least one of
the numbers $h_{1,2}^-$ is positive and at least one of the numbers $h_{1,2}^+$ is negative.

We assume that the system $u'=f(u)$ has a heteroclinic trajectory $\gamma$ crossing $\Sigma$ at $w=\gamma(0)\in\Sigma$.
Then
\begin{equation}
\label{gamma(xi)}
    \gamma(\xi)
  = \begin{cases}
       u^- + e^{\xi Q^-}(w-u^-),\qquad \xi \le 0,\\
       u^+ + e^{\xi Q^+}(w-u^+),\qquad \xi\ge 0.
    \end{cases}
\end{equation}

\begin{theo}\label{theo:simple}
The TWS corresponding to $\gamma$ is spectrally stable:
$$
  \spec \calL\subset \{\lambda\in\mC:\re\lambda<0\}\cup\{0\}.
$$
\end{theo}

Theorem \ref{theo:simple} holds for all types of shocks.
Next we give an example of stability loss for overcompressive shocks:
$u^-$ is a source and $u^+$ a sink.
Let $\Sigma = \{u_2=0\}$ and $u^\pm$ be as in (\ref{u+-}).

Suppose the matrices $Q^\pm$ have different eigendirections:
\begin{eqnarray}
\label{eq:Q-}
 &Q^- = \diag(h_1^-,h_2^-), \qquad h_j^\pm>0,\\
 \label{eq:Q+}
 &Q^+ =
 \left(
 \begin{array}{cc}  a_+ & c_+\\
                                 c_+ & b_+
              \end{array}
              \right)
              =T^{-1}
              \left(
\begin{array}{cc}
h_1^+ & 0\\
0 & h_2^+
\end{array}
\right)T,
\qquad
T=\left(
\begin{array}{cc}
\cos\theta & -\sin\theta\\
\sin\theta & \cos\theta
\end{array}
\right),
\end{eqnarray}
where  $0>h_2^+>h_1^+$ and $0<\theta<\pi/2$. Then
$$
a_+=h_1^+\cos^2\theta+h_2^+\sin^2\theta<0 ,\quad
b_+=h_2^+\cos^2\theta+h_1^+\sin^2\theta<0,\quad
c_+=(h_2^+-h_1^+)\sin\theta\cos\theta>0.
$$

The phase portrait of the system $u'=f(u)$
is presented on Figure \ref{fig:Over}.
A heteroclinic $\gamma=\gamma_s$ joining $u^-$ and $u^+$ and crossing $\Sigma$ at $w=(s,0)$ exists for
\begin{equation}
\label{eq:s0}
   s > s_0 = u_1^++\frac{b_+ u_2^+}{c_+}.
\end{equation}

As always for overcompressive shocks, 0 is a double eigenvalue of the spectral problem associated with $\gamma_s$.

If all parameters are fixed except $\theta$, and $\theta>0$ is small enough, then
$\gamma_s$ is spectrally stable by  Theorem~\ref{theo:simple}. Next we fix the matrices $Q^\pm$
and consider the loss of stability when $u^\pm$ and $s$ are changed.

\begin{theo}
\label{theo:over}
There exists a constant $\chi$ such that:
\begin{itemize}
\item
For $u_1^+-u_1^->\chi$,
all TWSs $\gamma^s$, $s>s_0$, are spectrally unstable: there exists an eigenvalue $\lambda>0$.

\item
Let $u_1^+-u_1^-<\chi$. Then there exists $s_2>s_1>s_0$ such that:
\begin{itemize}
\item
for $s_0<s<s_1$ the spectral problem corresponding to $\gamma_s$
has no real eigenvalues $\lambda>0$;
\item
for $s=s_1$ two eigenvalues $\lambda(s)$, $\bar\lambda(s)$
cross the imaginary axis: $\re\lambda(s)<0$ for $s_1-\delta<s<s_1$ and
$\re\lambda(s_1)=0$, $\im \lambda(s_1)>0$;
\item
for $s>s_2$ there exists  a real positive eigenvalue $\lambda= cs^{-1}+O(s^{-2})$,
$c>0$.
\end{itemize}
\end{itemize}
\end{theo}

Finally we give an example where instability appears as a result of a bifurcation.  Let
$$
u^- = \begin{pmatrix}0\\-2\end{pmatrix},\quad
  u^{\times\eps} = \begin{pmatrix}\eps\\\eps\end{pmatrix},\quad
  u^+ = \begin{pmatrix}-2\\\chi\end{pmatrix},
$$
and
$$
f^\eps
  =  \left\{ \begin{array}{lcl}
               Q^-(u-u^-)\quad                 & \mbox{ if } & \mbox{ $u_1 > -1$ and $u_2 < -1$}, \\
               Q^\times(u-u^{\times\eps})\quad & \mbox{ if } & \mbox{ $u_1 > -1$ and $u_2 > -1$}, \\
               Q^+(u-u^+)\quad                 & \mbox{ if } & \mbox{ $u_1 < -1$ },
            \end{array}
     \right.
     $$
where
\begin{eqnarray*}
 &    Q^- = \diag(\nu^-,\kappa^-), \quad
 Q^\times = \diag(\kappa^\times,\nu^\times), \quad
      Q^+ = \diag(\kappa^+,\nu^+), & \\
 &    \nu^- < 0 < \kappa^-, \quad
      \nu^\times < 0 < \kappa^\times, \quad
      \nu^+ < \kappa^+ < 0 . &
\end{eqnarray*}
The phase portraits of system (\ref{ode}) for $\eps=0$ and for
small $\eps > 0$  are presented on Figure \ref{fig:Bifur}.
The vector field $f^\eps$ is discontinuous on the lines
$$
  \Sigma^- = \{u_1 = -1\} \quad\mbox{and}\quad \Sigma^+ = \{u_2 = -1\}.
$$
The definition of solutions for system (\ref{ode}) does not meet any problem
since trajectories cross  the discontinuity lines transversely.
For $\eps = 0$ we define $\gamma^-$ as a heteroclinic from $u^-$ to $u^{\times\, 0}$
and $\gamma^+$ as a heteroclinic from $u^{\times\, 0}$ to $u^+$.
For small $\eps > 0$ there exists a heteroclinic $\gamma^\eps$ from $u^-$ to $u^+$.
For $\eps\to 0$ it converges pointwise to the concatenation of $\gamma^-$ and $\gamma^+$.
As in the previous section, we use smoothing to define  the spectral problem associated to $\gamma^\eps$.

\begin{theo}
\label{theo:main}
\begin{itemize}
\item
For $\eps=0$ the spectra of the  eigenvalue problems associated with the TWS's $\gamma^\pm$
lie in  $\{\lambda\in\mC : \re\lambda< 0\}\cup\{0\}$.
\item
For any small $\eps>0$ the eigenvalue problem associated with  $\gamma^\eps$
has an  eigenvalue
$$
  \lambda = c\eps + o(\eps), \qquad
  c = \frac{\nu^- \kappa^{\times\,2} (2\nu^\times + \chi)}{4\nu^\times (\kappa^\times - \nu^-)}.
$$
\end{itemize}
\end{theo}

The first statement is a particular case of Theorem \ref{theo:simple}.

\begin{cor}
If $\chi  > -2\nu^\times$, then $c>0$. This implies instability of  $\gamma^\eps$
for small $\eps>0$.
\end{cor}

Theorems  \ref{theo:simple} and \ref{theo:over} are proved in Section \ref{sec:simple},
and Theorem \ref{theo:main}  in Section \ref{sec:proof}.

\section{Shocks with a  single discontinuity}
\label{sec:simple}

In this section we consider stability of shocks generated by a heteroclinic
$\gamma$  which crosses the singular hyperplane $\Sigma$ just once.
Let $U^\pm\subset\mR^n$ be domains separated by  $\Sigma$, and
$$
  f^\pm = f|_{U^\pm}=Q^\pm (u-u^\pm),\qquad u^\pm \in U^\pm .
$$
We assume that the matrices $Q^\pm$ are symmetric.

A heteroclinic trajectory (\ref{gamma(xi)}) joining $u^-$ with $u^+$
is determined by a crossing point $w=\gamma(0)\in \Sigma$.
Let
$$
\hat f^\pm =f^\pm(w)=Q^\pm(w-u^\pm)
$$
and let $\bn$ be the unit normal to $\Sigma$ from $U^-$ to $U^+$. Then
$\gamma$ will cross $\Sigma$ transversely provided
$$
\langle \hat f^\pm,\bn\rangle >0.
$$

The eigenvalue  system (\ref{eq:L+-}) has the form
\begin{equation}
\label{eq:xv}
  v' = Q^\pm v+z,\quad  z' = \lambda v, \qquad
  \pm \xi > 0.
\end{equation}
For $\lambda=0$ we obtain $z\equiv z^0$, where $z^0=0$ because we need $z(\pm \infty)=0$.
Then
$$
  v(\xi)=e^{\xi Q^\pm}v^\pm,\qquad \pm \xi >0.
$$
We have $v(\pm\infty)=0$ if
$$
  v^+ = v(+0)\in E^s(Q^+), \quad
  v^- = v(-0)\in E^u(Q^-),
$$
where $E^{s,u}(Q^\pm)\subset\mR^n$ denotes the stable and unstable subspaces of a matrix:
$E^s$ is spanned by eigenvectors associated to negative eigenvalues, and $E^u$ -- to positive.

By Proposition \ref{prop:smooth}, $v^\pm$ satisfy (\ref{eq:S}). Hence $\lambda=0$ is an eigenvalue
with multiplicity
$$
  \dim E^s(Q^+)\cap S E^u(Q^-)\ge 1.
$$

\begin{rem}
If all eigenvalues of $Q^-$ are positive, and all eigenvalues of $Q^-$ are negative (fully overcompressive shock),
then for any nonzero $v^-\in\mR^n$, setting $v^+=Sv^-$ we get an eigenfunction.
Thus $\lambda=0$ is an eigenvalue of multiplicity $n$, as it should be.
\end{rem}

From now on we consider nonzero eigenvalues with $\re\lambda\ge 0$.

\begin{lem}\label{lem:P}
Let $\re\lambda\ge 0$. Then there exist unique matrices $P^\pm(\lambda)$, analytically depending on $\lambda$,
such that $\mp\re\sigma>0$ for all $\sigma\in \spec P^\pm$ and
\begin{equation}
\label{eq:Ppm}
(P^\pm)^2-P^\pm Q^\pm-\lambda I=0.
\end{equation}
\end{lem}

\noindent
{\it Proof}.
Since $Q^\pm$ are symmetric, there exist matrices  $T_\pm\in SO(n)$ such that
$$
  Q^\pm = T_\pm^{-1}\diag(h_1^\pm,\dots,h_n^\pm)T_\pm, \qquad
  h_{j}^\pm \ne 0.
$$
Then $P^\pm = T_\pm^{-1} \diag(\sigma_1^\pm,\dots,\sigma_n^\pm)T_\pm$, where
\begin{equation}
\label{eq:quad}
\sigma^2-h_j^\pm\sigma-\lambda=0.
\end{equation}
As proved in Section \ref{sec:quad}, if $\re\lambda\ge 0$, these
equations have unique solutions $\sigma_j^\pm$ with $\mp\re\sigma_j^\pm>0$,
and they analytically depend on $\lambda$ in the right half plane:
\begin{equation}
\label{eq:sigma}
\sigma_j^\pm=\frac12\Big(h_j^\pm\mp\sqrt{(h_j^\pm)^2+4\lambda}\Big).
\end{equation}
The value of the root must be taken with $\re\sqrt{\cdot}>0$.
\qed

\medskip

For large $|\lambda|$, the matrices $P^\pm$ are given by
\begin{equation}
\label{eq:infty}
P^\pm(\lambda)=\mp\sqrt{\lambda}I+\frac12 Q^\pm+O(|\lambda|^{-1/2}),
\end{equation}
where $\sqrt{\lambda}$ is the branch of root with $\re\sqrt{\lambda}\ge 0$ for $\re\lambda\ge 0$.

If the shock is fully overcompressive, i.e.\ $\mp h_j^\pm>0$, then for small $\lambda$,
\begin{equation}
\label{eq:lambda0}
  P^\pm(\lambda) = Q^\pm+\lambda (Q^\pm)^{-1}+O(|\lambda|^2) .
\end{equation}

\begin{lem}
For $\re\lambda\ge 0$, $\lambda \ne 0$, solutions of the eigenvalue problem (\ref{eq:xv}) have the form
\begin{eqnarray*}
v(\xi)=e^{\xi P^-}v^-,\quad  z(\xi)=e^{\xi P^-}z^0,\qquad \xi<0,\\
v(\xi)=e^{\xi P^+}v^+,\quad  z(\xi)=e^{\xi P^+}z^0,\qquad \xi>0,
\end{eqnarray*}
where
\begin{equation}
\label{eq:Pz}
   P^\pm z^0 = \lambda v^\pm, \quad v^+=Sv^-.
\end{equation}
\end{lem}

\noindent
{\it Proof}. As in the proof of Lemma \ref{lem:P},
we may  assume that the matrices $Q^\pm$ are diagonal.
Then for the components of  $v(\xi)$ we obtain the equations
$$
 v''_j-h_j^\pm v'_j-\lambda v_j=0,\qquad \mp\xi>0.
$$
Since we need a solution with $v^\pm(\pm\infty)=0$,
we must take the root $\sigma=\sigma_j^\pm$ of equation (\ref{eq:quad}) with $\mp\re\sigma>0$. Then
$$
v_j(\xi)=e^{\xi \sigma^\pm_j}v_j^\pm,\quad z_j(\xi)=\frac{\lambda e^{\xi \sigma^\pm_j}v_j^\pm}{\sigma_j^\pm},
\qquad \mp\xi>0.
$$
Thus $v_j^\pm=\lambda^{-1}\sigma^\pm_j z_j^0$.
\qed

\medskip

By (\ref{eq:Pz}), the spaces (\ref{eq:Epm}) are given by
$$
     E^\pm(\lambda)
  =  \Big\{ \bv = \begin{pmatrix} v^\pm \\ z^0 \end{pmatrix}:
                         z^0 = \lambda {P^\pm}^{-1} v^\pm \Big\}.
$$

If  $\lambda\ne 0$, $\re\lambda\ge 0$, is an eigenvalue of the spectral problem,
then by (\ref{eq:Pz}) the corresponding $z^0\ne 0$ must satisfy the equation
$$
  P^+z^0 = SP^-z^0.
$$
We obtain

\begin{prop}
$\lambda\ne 0$, $\re\lambda\ge 0$, is an eigenvalue  iff
\begin{equation}
\label{eq:Delta}
D(\lambda)=\det(P^+(\lambda)-SP^-(\lambda))=0.
\end{equation}
\end{prop}

\begin{rem}
$D(\lambda)$ is an analog of the Evans function $\calD(\lambda)$  \cite{Evans} for piece-wise linear vector fields.
However, they are slightly different: $\calD(\lambda)$ has a zero at $\lambda=0$,
while in $D(\lambda)$ this degeneracy is removed.
\end{rem}

\begin{rem}
For given $Q^\pm$ and $S$  equations (\ref{eq:Ppm}) and (\ref{eq:Delta})
can be rewritten as a
system of  polynomial  equations in $2n+1$ variables  $\sigma_{j}^\pm$, $\lambda$.
The number of complex solutions is given by Bezout's theorem.
However conditions $\mp\re\sigma_j^\pm>0$, $\re\lambda\ge 0$,
decrease the number of possible solutions.
\end{rem}

\medskip

Equation (\ref{eq:infty}) implies

\begin{lem}
\label{lem:infty}
As $|\lambda|\to +\infty$, $\re\lambda\ge 0$, we have
$$
 D(\lambda)= (-\sqrt{\lambda})^n\big(c+O(|\lambda|^{-1/2})\big),\qquad c=\det(I+S)>0.
$$
\end{lem}

\label{sec:proof12}

\subsection{Proof of Theorem \ref{theo:simple}}
\label{sec:diag}

From now on we assume $u\in\mR^2$, then  $\Sigma$ is a line.
Rotating the coordinate system, we may assume $\Sigma=\{u_2=0\}$ and
$\bn =(0,1)$. Then the matrix $S$ in (\ref{eq:S}) is given by
\begin{equation}
\label{eq:alpha-beta}
  S =
  \begin{pmatrix}
  1 & \alpha \\
  0 & \beta
              \end{pmatrix}, \qquad
  \alpha = \frac{\hat f_1^+ - \hat f_1^-}{\hat f_2^-}, \quad
   \beta = \frac{\hat f_2^+}{\hat f_2^-}.
\end{equation}

Consider the simple case when $\bn$ is an eigenvector for both $Q^\pm$, then
\begin{equation}
\label{eq:P-}
    Q^\pm
 = \diag(  h_1^\pm,h_2^\pm),\quad
    P^\pm
 = \diag(\sigma_1^\pm,\sigma_2^\pm),\qquad
   \mp \re \sigma_{j}^\pm(\lambda) > 0.
\end{equation}
For overcompressive shocks we have $h_j^->0$ and $h_j^+<0$,
while for Lax shocks $h_1^-h_2^-<0$ and $h_j^+<0$.
Overcompressive shocks appear in families, see Figure \ref{fig:Over},
while for Lax shocks there is a unique heteroclinic.

The next proposition implies Theorem \ref{theo:simple}.

\begin{prop}\label{prop:stable}
If $\bn$ is an eigenvector for both $Q^\pm$, then for all heteroclinics joining $u^-$ and $u^+$,
there  are no eigenvalues of the spectral problem with $\re\lambda\ge 0$ except  $\lambda=0$.
The shock is spectrally stable.
\end{prop}

\noindent
{\it Proof.}  By (\ref{eq:Delta}),
\begin{equation}
\label{eq:D}
 D(\lambda)=\left|
\begin{array}{cc}
\sigma_1^+-\sigma_1^- & -\alpha \sigma_2^-\\
0 & \sigma_2^+-\beta \sigma_2^-
\end{array}
\right|
=(\sigma_1^+-\sigma_1^-)(\sigma_2^+-\beta \sigma_2^-).
\end{equation}
Since $\beta>0$, we have
$$
\re (\sigma_1^+-\sigma_1^-)<0,\quad \re(\sigma_2^+-\beta \sigma_2^-)<0.
$$
Hence $ D(\lambda)\ne 0$. For Lax or undercompressive shocks $\lambda=0$ is a simple eigenvalue,
and for overcompressive shocks a double one.
\qed

\subsection{Overcompressive shocks}
\label{sec:over}

In this section we consider overcompressive shocks. Then
$$
  Q^\pm = \begin{pmatrix}  a_\pm & c_\pm\\
                                   c_\pm & b_\pm
                \end{pmatrix},
$$
where
$$
  \mp a_\pm>0,\quad \mp b_\pm>0,\quad \det Q^\pm=a_\pm b_\pm-c_\pm^2>0.
$$
Heteroclinics $\gamma=\gamma_s$ appear in families determined by the  break point  $\gamma_s(0)=w=(s,0)\in\Sigma$.
Then
$$
  \hat f^\pm=Q^\pm w -Q^\pm u^\pm
$$
has components
$$
  \hat f_2^\pm = (s-u_1^\pm) c_\pm-b_\pm u_2^\pm,\quad
  \hat f_1^\pm = (s-u_1^\pm) a_\pm - c_\pm u_2^\pm.
$$
Admissible values of $s$ are such that $\hat f_2^\pm>0$. Then
\begin{eqnarray}
\label{eq:beta}
      \hat f_2^-\beta
 &=&  s c_+ - b_+u_2^+-c_+u_1^+, \\
      \hat f_2^-\alpha
 &=&  s(a_+-a_-) + a_-u_1^- - a_+u_1^+ + c_-u_2^- - c_+u_2^+ .
\label{eq:alpha}
\end{eqnarray}

\begin{lem}
\label{lem:Delta0}
\begin{equation}
\label{eq:Delta(0)}
\hat f_2^-D(0)=\Lambda(u^-,u^+):=A_1(u_1^+-u_1^-)+A_2(u_2^+-u_2^-),
\end{equation}
where
\begin{equation}
\label{eq:A}
A_1=c_-\det Q^+-c_+\det Q^-,\quad A_2=b_{-}\det Q^+-b_+ \det Q^-.
\end{equation}
Thus $A_2>0$, while $A_1$ may have any sign.
\end{lem}

The proof is a direct computation:
\begin{eqnarray*}
D(0)&=&\det(Q^+-SQ^-)=
\left|
\begin{array}{cc}
a_+-a_- - \alpha c_- & c_+-c_- -\alpha b_-\\
c_+-\beta c_- & b_+-\beta b_-
\end{array}
\right|\\
&=&\alpha(c_+b_--c_-b_+)+\beta \big( c_-(c_+-c_-)-b_-(a_+-a_-)\big)+(a_+-a_-)b_+-(c_+-c_-)c_+.
\end{eqnarray*}
Substituting $\alpha$ and $\beta$ from (\ref{eq:beta})--(\ref{eq:alpha}), we see that the terms
in $\hat f_2^-D(0)$ with $s$ cancel, and we obtain (\ref{eq:Delta(0)}).
\qed

\medskip

Note that $\hat f_2^-D(0)$  has the same value
for all heteroclinics joining $u^-$ and $u^+$.
We obtain

 \begin{equation}
\label{eq:FG}
\hat f_2^- D(\lambda)=sF(\lambda)+G(\lambda),\qquad F(0)=0,\quad G(0)=\Lambda(u^-,u^+).
\end{equation}

\begin{cor}\label{cor:unstable}
If $\Lambda(u^-,u^+) < 0$, then for every heteroclinic $\gamma_s$ joining $u^-$ and $u^+$,
the corresponding spectral problem has an eigenvalue $\lambda>0$,
and hence the shock is spectrally unstable.
\end{cor}

Indeed, then $D(0)<0$ and $D(+\infty)=+\infty$ by (\ref{eq:infty}).

\begin{cor}
If $\Lambda(u^-,u^+) > 0$,\footnote{This holds, for example, if $u_1^+=u_1^-$.}
loss of spectral stability with real $\lambda$ crossing 0 is impossible.
As $s$ changes, loss of stability of $\gamma_s$ may happen only when
a pair of complex conjugate eigenvalues crosses the imaginary axis.
\end{cor}

\subsection{Proof of Theorem \ref{theo:over}}

Consider the simplest nontrivial case:
$\bn$ is an eigenvector for $Q^-$ but not for $Q^+$.
Let $Q^\pm$ be as in (\ref{eq:Q-})--(\ref{eq:Q+}).
Then we have $\hat f_2^-=-h_2^-u_2^->0$ and
$$
  \hat f_2^+=(s-u_1^+)c_+ - b_+u_2^+,\qquad c_+>0,
$$
is positive for $s$ as in (\ref{eq:s0}).
Only for such values of $s$  there exists a heteroclinic $\gamma_s$  passing through $w=(s,0)$,
see Figure \ref{fig:Over}.

In (\ref{eq:A}) we have
$$
  A_1=-c_+ h_1^-h_2^-<0,\qquad A_2=h_2^-h_1^+h_2^+-b_+h_1^-h_2^->0.
$$
Hence for
$$
  u_1^+-u_1^- >\chi=\frac{h_2^-h_1^+h_2^+-b_+h_1^-h_2^-}{c_+ h_1^-h_2^-},
$$
we have $\Lambda(u^-,u^+)<0$, and then by Corollary \ref{cor:unstable} all heteroclinics $\gamma_s$
are unstable shocks. This proves the first statement of Theorem \ref{theo:over}.
For $ u_1^+-u_1^- <\chi$, we have $\Lambda(u^-,u^+)<0$.

Let $\sigma_j^\pm$ be as in (\ref{eq:sigma}). Then
$$
    P^-
  = \begin{pmatrix} \sigma_1^- & 0\\
                                     0 & \sigma_2^-
          \end{pmatrix},  \quad
    P^+
  = \begin{pmatrix} A_+ & C_+ \\
                             C_+& B_+
          \end{pmatrix},
$$
where
$$
  A_+ = \sigma_1^+\cos^2\theta+\sigma_2^+\sin^2\theta,\quad
  B_+ = \sigma_2^+\cos^2\theta+\sigma_1^+\sin^2\theta,\quad
  C_+ = (\sigma_2^+-\sigma_1^+)\sin\theta\cos\theta.
$$

We obtain
\begin{eqnarray*}
      D(\lambda)
  &=& \det(P^+-SP^-)=
      \left|\begin{array}{cc}
             A_+-\sigma_1^- & C_+-\alpha \sigma_2^-\\
                        C_+ & B_+-\beta \sigma_2^-
            \end{array}\right| \\
  &=& \alpha C_+\sigma_2^- - \beta \sigma_2^-(A_+-\sigma_1^-)+(A_+-\sigma_1^-)B_+ - C_+^2.
\end{eqnarray*}
By (\ref{eq:beta})--(\ref{eq:alpha}), the function $F$ in (\ref{eq:FG}) is given by
\begin{eqnarray}
F&=&(a_+-h_1^-) C_+ \sigma_2^- -  c_+ (A_+-\sigma_1^-)\sigma_2^- \nonumber \\
&=&\sigma_2^-\sin\theta\cos\theta\Big((h_1^+\cos^2\theta+h_2^+\sin^2\theta-h_1^-)(\sigma_2^+-\sigma_1^+)\nonumber\\
&&\qquad -(h_2^+-h_1^+)(\sigma_1^+\cos^2\theta+\sigma_2^+\sin^2\theta-\sigma_1^-)\Big).
\label{eq:F}
\end{eqnarray}

For small $\lambda$ we  have
$$
\sigma_j^\pm=h_j^\pm+\lambda/h_j^\pm+O(|\lambda|^2).
$$
Substituting in (\ref{eq:F}) and performing  simplifications, we obtain
\begin{eqnarray}
F'(0)=\frac{(h_2^+-h_1^+)(h_1^--h_2^+)(h_1^--h_1^+)\sin\theta\cos\theta}{h_1^+h_1^-}<0.
\end{eqnarray}

Let $0<\eps=s^{-1}\ll 1$. Then the equation
$$
sF(\lambda)+G(\lambda) = 0\quad\iff\quad F(\lambda)+\eps G(\lambda) = 0
$$
has a solution
$$
  \lambda = -\eps\frac{G(0)}{F'(0)}+O(\eps^2)>0.
$$
We obtain

\begin{prop}
If $u_1^+-u_1^-<\chi$, then for sufficiently large $s>s_0$ the shock  is spectrally unstable:
there exists a real eigenvalue $\lambda=c s^{-1}+O(s^{-2})$, $c>0$.
\end{prop}

\begin{cor}
Suppose $h_j^\pm$ and $u^\pm$ are fixed while $\theta>0$ is small enough.
Then there exist $s_0<s_1<s_2$ such that for $s_0<s<s_1$ the shock is stable.
For $s=s_1$ a pair of complex eigenvalues $\lambda_\pm=\mu(s)\pm i \omega(s)$ crosses
the imaginary axis: $\mu(s)<0$ for $s_1-\delta<s<s_1$, and $\mu(s_1)=0$, $\omega(s_1)>0$.
For $s>s_2$ there exists a real eigenvalue $\lambda= c s^{-1}+O(s^{-2})$, $c>0$.
\end{cor}

Now Theorem \ref{theo:over} is proved.

\section{Proof of Theorem \ref{theo:main}}
\label{sec:proof}

\subsection{The curves $\gamma^\pm$ and $\gamma^\eps$}

Equations (\ref{ode}) can be solved explicitly. In particular, we obtain:
$$
    \gamma^-(\xi)
  = \left\{ \begin{array}{cl}
            \bigg(\!\!\! \begin{array}{c} 0 \\ -2 + e^{\kappa^-\xi}
                         \end{array} \!\!\!
            \bigg) \quad & \mbox{if } \xi < 0, \\
            \bigg(\!\!\! \begin{array}{c} 0 \\ - e^{\nu^\times\xi}
                         \end{array} \!\!\!
            \bigg) \quad & \mbox{if } \xi > 0,
            \end{array}
    \right. \quad
    \gamma^+(\xi)
  = \left\{ \begin{array}{cl}
            \bigg(\!\!\! \begin{array}{c} - e^{\kappa^\times \xi} \\ 0
                         \end{array} \!\!\!
            \bigg) \quad & \mbox{if } \xi < 0, \\
            \bigg(\!\!\! \begin{array}{c} -2 + e^{\kappa^+\xi} \\ \chi - \chi e^{\nu^+\xi}
                         \end{array} \!\!\!
            \bigg) \quad & \mbox{if } \xi > 0.
            \end{array}
    \right.
$$
To present explicit expressions for $\gamma^\eps$, we introduce some notation. We define $\eps^\times$ and two constants $\xi^- < \xi^+$ (in fact, we fix only the difference $\xi^+ - \xi^-$) such that
\begin{equation}
\label{epstimes}
  \eps^\times = \eps - (1+\eps) \Big(\frac{\eps}{1+\eps}\Big)^{-\nu^\times/\kappa^\times}, \quad
  e^{\kappa^\times (\xi^+ - \xi^-)} = \frac{1+\eps}{\eps}
\end{equation}
Then by a direct computation we obtain:
$$
    \gamma^\eps(\xi)
  = \left\{ \begin{array}{cl}
            \bigg(\!\!\! \begin{array}{c} 0 \\ -2 + e^{\kappa^-(\xi-\xi^-)}
                         \end{array} \!\!\!
            \bigg) \quad & \mbox{if } \xi \le \xi^-, \\
            \bigg(\!\!\! \begin{array}{c} \eps - \eps e^{\kappa^\times(\xi-\xi^-)} \\
                                         \eps - (1+\eps) e^{\nu^\times(\xi-\xi^-)}
                         \end{array} \!\!\!
            \bigg) \quad & \mbox{if } \xi^- \le \xi \le \xi^+, \\
            \bigg(\!\!\! \begin{array}{c} - 2 + e^{\kappa^+(\xi-\xi^+)} \\
                                        \chi - (\chi - \eps^\times) e^{\nu^+(\xi-\xi^+)}
                         \end{array} \!\!\!
            \bigg) \quad & \mbox{if } \xi^+ \le \xi.
            \end{array}
    \right.
$$

The curve $\gamma^\eps$ is piece-wise smooth: ${\gamma^\eps}'(\xi)$ is discontinuous at $\xi=\xi^\pm$.

\subsection{Operator $\bL$}

We put $\bv = \begin{pmatrix} v(\xi)\\ z(\xi) \end{pmatrix}$.
Then  the spectral problem (\ref{eq:system}) takes the form $\bL\bv = 0$, where for $\xi\ne\xi^\pm$,
$$
    \bL
  = \partial_\xi - \bQ_\lambda^*, \qquad
    \bQ_\lambda^*
  = \bigg(\!\begin{array}{cc} Q^* & I \\ \lambda I & 0 \end{array}\!\bigg),
$$
where $*$ equals $-, \times$, or $+$. The passage through $\xi^\pm$ is performed with the help of Proposition \ref{prop:smooth}:
$$
    \bL
  = \left\{\begin{array}{lcc}
        \partial_\xi - \bQ^-_\lambda            &\mbox{if }& \xi < \xi^- - 0, \\[1mm]
        \bv(\xi^- + 0)
        = \bS^- \bv(\xi^- - 0),  && \\[1mm]
        \partial_\xi - \bQ^\times_\lambda       &\mbox{if }& \xi^- + 0 < \xi < \xi^+ - 0, \\[1mm]
        \bS^+ \bv (\xi^+ + 0)
        = \bv(\xi^+ - 0),  && \\[1mm]
        \partial_\xi - \bQ^+_\lambda            &\mbox{if }& \xi^+ + 0 < \xi,
          \end{array}\right.
$$
where
\begin{eqnarray*}
&   \bS^\pm
  = \left( \begin{array}{cc} S^\pm & 0 \\  0 & I
           \end{array}\right), \quad
    S^-
  = \left( \begin{array}{cc} 1 & -\eps\kappa^\times / \kappa^- \\
                             0 & -(1+\eps)\nu^\times / \kappa^-
           \end{array}\right), \quad
    S^+
  = \left( \begin{array}{cc}
                  - (1+\eps)\kappa^\times / \kappa^+    & 0 \\
      (\chi\nu^+ + \eps^\flat\nu^\times) / \kappa^+     & 1
           \end{array}\right). &
\end{eqnarray*}
The number $\eps^\flat$ is defined by
$$
  \eps^\flat\nu^\times  =  - \eps^\times \nu^+ + (\eps^\times - \eps) \nu^\times.
$$

The $2\times 2$ matrices $S^\pm$ should be computed with the help of (\ref{eq:S-})--(\ref{eq:S+}).
It is important to note that the operator $\bS^+$ maps $\bv(\xi^+ + 0)$ to $\bv(\xi^+ - 0)$, not vise versa.
Because of this $S^+$ is the  inverse of the matrix determined by (\ref{eq:S+}).

First, we compute $S^-$. Then $\Sigma=\{u_2=-1\}$ and $\bn=(0,1)$. Then
\begin{eqnarray*}
      \hat f^-
  &=& Q^- \bigg(\bigg(\!\!\begin{array}{c} 0 \\ -1 \end{array}\!\!\bigg) - u^- \bigg)
   =  \bigg(\!\!\begin{array}{c} 0 \\ \kappa^- \end{array}\!\!\bigg), \\
      \Delta f
  &=& Q^\times \bigg(\bigg(\!\!\begin{array}{c} 0 \\ -1 \end{array}\!\!\bigg) - u^{\times\eps} \bigg)
    - Q^- \bigg(\bigg(\!\!\begin{array}{c} 0 \\ -1 \end{array}\!\!\bigg) - u^- \bigg)
   =  - \bigg(\!\!\begin{array}{c} \eps\kappa^\times \\ (1 + \eps)\nu^\times - \kappa^- \end{array}\!\!\bigg).
\end{eqnarray*}
By using (\ref{eq:S-}) we obtain $S^-$.

To compute $S^+$, we take $\bn = (-1,0)$. Then
\begin{eqnarray*}
     \hat f^-
 &=& Q^\times \bigg(\bigg(\!\!\begin{array}{c} -1 \\ \eps^\times \end{array}\!\!\bigg) - u^{\times\eps} \bigg)
  =  \bigg(\!\!\begin{array}{c} -(1+\eps)\kappa^\times \\ (\eps^\times - \eps)\nu^\times \end{array}
     \!\!\bigg), \\
      \Delta f
  &=& Q^+ \bigg(\bigg(\!\!\begin{array}{c} -1 \\ \eps^\times \end{array}\!\!\bigg) - u^+ \bigg)
      - \hat f^-
   =  \bigg(\!\!\begin{array}{c} \kappa^+ + (1+\eps)\kappa^\times \\
                               (\eps^\times - \chi)\nu^+ - (\eps^\times - \eps)\nu^\times \end{array}\!\!\bigg).
\end{eqnarray*}
We have: $S^+ = S^{-1}$, where $S$ satisfies (\ref{eq:S+}):
$$
    S^+
  = \bigg(\!\!\begin{array}{cc} -\kappa^+ / ((1 + \eps) \kappa^\times)  &  0  \\
      (\chi\nu^+ + \eps^\flat\nu^\times) / ((1 + \eps) \kappa^\times)   &  1
              \end{array}\!\!\bigg)^{-1}
  = \left( \begin{array}{cc}
                  - (1+\eps)\kappa^\times / \kappa^+    & 0 \\
      (\chi\nu^+ + \eps^\flat\nu^\times) / \kappa^+     & 1
           \end{array}\right).
$$

\subsection{Eigenfunction of $\bL$}

Let $\bv = \begin{pmatrix} v(\xi) \\ z(\xi)\end{pmatrix}$ be an eigenfunction of $\bL$
with eigenvalue $\lambda$, $\re\lambda > 0$. We assume that $|\lambda|$ is small.

\subsubsection{The function $\bv$ for $\xi < \xi^-$ and $\xi > \xi^+$}

The ODE system $\partial_\xi\bv - \bQ_\lambda^-\bv = 0$ splits into two linear systems with constant coefficients:
$$
  v_1' = \nu^- v_1 + z_1, \quad z_1' = \lambda v_1, \qquad
  v_2' = \kappa^- v_2 + z_2, \quad z_2' = \lambda v_2.
$$
The corresponding characteristic equations are
\begin{equation}
\label{quad}
  \sigma_1^2 - \nu^- \sigma_1 - \lambda = 0, \quad
  \sigma_2^2 - \kappa^- \sigma_2 - \lambda = 0.
\end{equation}
The condition $\bv(\pm\infty)=0$ imposes the restrictions $\re\sigma_1 > 0$ and $\re\sigma_2 > 0$. We have the following obvious lemma (see also Section \ref{sec:quad})

\begin{lem}
\label{lem:-obv}
If $\re\lambda>0$ and $|\lambda|$ is small then any equation (\ref{quad}) has a unique solution with positive real part:
\begin{equation}
\label{sigma-}
   \sigma_1^- = -\lambda / \nu^- + O(\lambda^2), \quad
   \sigma_2^- = \kappa^- + \lambda / \kappa^- + O(\lambda^2).
\end{equation}
\end{lem}

Therefore for $\xi < \xi^-$,
$$
    v(\xi)
  = \bigg(\!\!\! \begin{array}{c}
            a_1 e^{\sigma_1^- (\xi-\xi^-)} \\
            a_2 e^{\sigma_2^- (\xi-\xi^-)}
                 \end{array}\!\!\!\bigg), \quad
    z(\xi)
  = \bigg(\!\!\! \begin{array}{c}
            a_1 (\sigma_1^- - \nu^-) e^{\sigma_1^- (\xi-\xi^-)} \\
            a_2 (\sigma_2^- - \kappa^-) e^{\sigma_2^- (\xi-\xi^-)}
                 \end{array}\!\!\!\bigg),
$$
where $a_1$ and $a_2$ are arbitrary constants.

Analogously the system $\partial_\xi\bv - \bQ_\lambda^+\bv = 0$ splits into two linear systems
$$
  v_1' = \kappa^+ v_1 + z_1, \quad z_1' = \lambda v_1, \qquad
  v_2' = \nu^+ v_2 + z_2, \quad z_2' = \lambda v_2.
$$
The corresponding characteristic equations are
\begin{equation}
\label{quad2}
  \sigma_1^2 - \kappa^+ \sigma_1 - \lambda = 0, \quad
  \sigma_2^2 - \nu^+ \sigma_2 - \lambda = 0.
\end{equation}
The condition $\bv(\pm\infty)=0$ implies $\re\sigma_1 < 0$ and $\re\sigma_2 < 0$.

\begin{lem}
\label{lem:+obv}
If $\re\lambda>0$ and $|\lambda|$ is small, then any equation (\ref{quad2}) has a unique solution with negative real part:
\begin{equation}
\label{sigma+}
  \sigma_1^+ = \kappa^+ + \lambda / \kappa^+ + O(\lambda^2), \quad
  \sigma_2^+ = \nu^+ + \lambda / \nu^+ + O(\lambda^2).
\end{equation}
\end{lem}

Therefore for $\xi > \xi^+$,
$$
    v(\xi)
  = \begin{pmatrix}
            b_1 e^{\sigma_1^+ (\xi-\xi^+)} \\
            b_2 e^{\sigma_2^+ (\xi-\xi^+)}
                 \end{pmatrix} , \quad
    z(\xi)
  =  \begin{pmatrix}
            b_1 (\sigma^+_1 - \kappa^+) e^{\sigma_1^+ (\xi-\xi^+)} \\
            b_2 (\sigma^+_2 - \nu^+) e^{\sigma_2^+ (\xi-\xi^+)}
                 \end{pmatrix},
$$
where $b_1$ and $b_2$ are arbitrary constants.

\begin{cor}
\label{cor:+-0}
\begin{eqnarray*}
& v(\xi^- - 0)  = a, \quad
  z(\xi^- - 0) = (P^- - Q^-) a, \quad
  v(\xi^+ + 0)  = b, \quad
  z(\xi^+ + 0) = (P^+ - Q^+) b, & \\
& a = \begin{pmatrix}
a_1 \\ a_2
\end{pmatrix}, \quad
  b =
  \begin{pmatrix}
  b_1 \\ b_2
  \end{pmatrix} ,
  \quad
      P^\pm
    = \begin{pmatrix}
    \sigma^\pm_1 & 0 \\
    0 & \sigma^\pm_2
    \end{pmatrix}. &
\end{eqnarray*}
\end{cor}

\subsubsection{The function $\bv$ for $\xi = \xi^- + 0$ and $\xi = \xi^+ - 0$}

We put $A = \bv(\xi^- + 0)$, $B = \bv(\xi^+ - 0)$. Then
$$
  A = \bS^- \left(\!\!\begin{array}{c} a \\ (\sigma^- - Q^-) a \end{array}\!\!\right)
    = \left(\!\!\begin{array}{c} S^- \\ P^- - Q^- \end{array}\!\!\right) a, \quad
  B = \bS^+ \left(\!\!\begin{array}{c} b \\ (P^+ - Q^+) b \end{array}\!\!\right)
    = \left(\!\!\begin{array}{c} S^+ \\ P^+ - Q^+ \end{array}\!\!\right) b.
$$

\subsubsection{The function $\bv$ for $\xi^- < \xi < \xi^+$}

In the interval $(\xi^-,\xi^+)$ the function $\bv$ satisfies the equation
$\partial_\xi\bv - \bQ_\lambda^\times \bv = 0$. This equation splits into two systems
$$
  v_1' = \kappa^\times v_1 + z_1, \quad z_1' = \lambda v_1, \qquad
  v_2' = \nu^\times v_2 + z_2, \quad z_2' = \lambda v_2.
$$
Characteristic polynomials are
$$
  {\sigma_1^\times}^2 - \kappa^\times\sigma_1^\times - \lambda = 0, \quad
  {\sigma_2^\times}^2 - \nu^\times\sigma_2^\times - \lambda = 0.
$$
For small $|\lambda|$ the roots are
\begin{equation}
\label{sigmatimes}
    \sigma_{11}^\times
  = -\frac{\lambda}{\kappa^\times} + O(\lambda^2), \quad
    \sigma_{12}^\times
  = \kappa^\times + O(\lambda), \quad
    \sigma_{21}^\times
  = -\frac{\lambda}{\nu^\times} + O(\lambda^2), \quad
    \sigma_{22}^\times
  = \nu^\times + O(\lambda).
\end{equation}
Consider the linear transformation $T$ which diagonalizes $\bQ_\lambda^\times$:
$$
    T \bQ_\lambda^\times T^{-1}
  = \calD_\lambda
 := \diag(\sigma_{11}^\times,\sigma_{21}^\times,\sigma_{12}^\times,\sigma_{22}^\times), \qquad
    T
  = \left(\begin{array}{cccc} \sigma_{11}^\times &      0              &   1  &   0        \\
                                    0            & \sigma_{21}^\times  &   0  &   1        \\
                             \sigma_{12}^\times  &      0              &   1  &   0        \\
                                    0            & \sigma_{22}^\times  &   0  &   1
          \end{array}\right) .
$$

The equation $\bL\bv = 0$ is equivalent to $B = e^{(\xi^+ - \xi^-) \bQ_\lambda^\times} A$. In more detail,
$$
    \left(\begin{array}{c} S^+ \\ P^+ - Q^+ \end{array} \right) b
  = e^{(\xi^+ - \xi^-) \bQ_\lambda^\times} \left(\begin{array}{c} S^- \\ P^- - Q^- \end{array} \right) a.
$$
This equation can be written as
$\Theta_0 \Big(\!\!\begin{array}{c} a \\ b \end{array} \!\!\Big) = 0$, where the $4\times 4$ matrix $\Theta_0$
has the form
$$
  \left(
        e^{(\xi^+-\xi^-) \bQ_\lambda^\times}
         \left( \begin{array}{c} S^- \\ P^- - Q^- \end{array} \right),
        - \left( \begin{array}{c} S^+ \\ P^+ - Q^+ \end{array} \right)
  \right).
$$

The equation $\bQ_\lambda^\times = T^{-1} \calD_\lambda T$ implies
$\Theta\left(\!\begin{array}{c} a\\ b \end{array}\!\right) = 0$,
where
$$
   \Theta
 = T\Theta_0
 = \left(
        e^{(\xi^+ - \xi^-) \calD_\lambda} T
         \left( \begin{array}{c} S^- \\ P^- - Q^- \end{array} \right),
        - T \left( \begin{array}{c} S^+ \\ P^+ - Q^+ \end{array} \right)
   \right).
$$

\subsubsection{Computation of $\Theta$}

We put
$$
  E = e^{(\xi^+-\xi_-)\calD_\lambda} = \diag(E_{11},E_{21},E_{12},E_{22}).
$$
Then by (\ref{epstimes}) and (\ref{sigmatimes})
\begin{eqnarray*}
     E_{11}
 &=& e^{(\xi^+ - \xi^-)\sigma_{11}^\times}
  =  e^{O(\lambda |\log\eps|)}, \\
     E_{21}
 &=& e^{(\xi^+ - \xi^-)\sigma_{21}^\times}
  =  e^{O(\lambda |\log\eps|)}, \\
     E_{12}
 &=& e^{(\xi^+ - \xi^-)\sigma_{12}^\times}
  =  \frac{1+\eps}{\eps} \big(1 +  O(\lambda |\log\eps|)\big), \\
      E_{22}
 &=& e^{(\xi^+ - \xi^-)\sigma_{22}^\times}
  =  \Big(\frac{1+\eps}{\eps}\Big)^{\nu^\times / \kappa^\times} \big(1 +  O(\lambda |\log\eps|)\big)
  =  \frac{\eps-\eps^\times}{1+\eps} \big(1 +  O(\lambda |\log\eps|)\big).
\end{eqnarray*}
Hence $E_{11}$ and $E_{21}$ are close to $1$, $E_{12}$ is approximately $1/\eps$, and $E_{22}$ is small.

We have
\begin{eqnarray*}
T \begin{pmatrix} S^- \\ P^- - Q^- \end{pmatrix}
 &=& \begin{pmatrix} \sigma_{11}^\times &           0         &      1    &      0        \\
                                    0            & \sigma_{21}^\times  &      0    &      1        \\
                              \sigma_{12}^\times &           0         &      1    &      0        \\
                                    0            & \sigma_{22}^\times  &      0    &      1
           \end{pmatrix}
     \begin{pmatrix}     1        &   -\eps\kappa^\times / \kappa^-         \\
                                 0        &   -(1+\eps)\nu^\times / \kappa^-        \\
                      \sigma_1^- - \nu^-  &                0                        \\
                                 0        &   \sigma_2^- - \kappa^-
           \end{pmatrix} \\
 &=& \begin{pmatrix}
 \sigma_{11}^\times + \sigma_1^- - \nu^- &  -\eps\sigma_{11}^\times \kappa^\times / \kappa^-   \\
                                 0       &  -(1+\eps) \sigma_{21}^\times \nu^\times / \kappa^-
                                                                      + \sigma_2^- - \kappa^-  \\
 \sigma_{12}^\times + \sigma_1^- - \nu^- &  -\eps\sigma_{12}^\times \kappa^\times / \kappa^-   \\
                                 0       &  -(1+\eps) \sigma_{22}^\times \nu^\times / \kappa^-
                                                                      + \sigma_2^- - \kappa^-
           \end{pmatrix}.
\end{eqnarray*}

and

\begin{eqnarray*}
T \begin{pmatrix} S^+ \\ P^+ - Q^+ \end{pmatrix}
 &=& \begin{pmatrix} \sigma_{11}^\times &          0         &      1    &      0    \\
                                        0            & \sigma_{21}^\times &      0    &      1    \\
                                  \sigma_{12}^\times &          0         &      1    &      0    \\
                                        0            & \sigma_{22}^\times &      0    &      1
              \end{pmatrix}
     \begin{pmatrix} -(1+\eps)\kappa^\times/\kappa^+  &            0          \\
                 (\chi\nu^+ + \eps^\flat\nu^\times) / \kappa^+    &            1          \\
                                \sigma_1^+ - \kappa^+             &            0          \\
                                            0                     &  \sigma_2^+ - \nu^+
           \end{pmatrix}
    \\
 &=& \begin{pmatrix}
  -(1+\eps) \sigma_{11}^\times \kappa^\times / \kappa^+  + \sigma_1^+ - \kappa^+ &        0              \\
  \sigma_{21}^\times (\chi\nu^+ + \eps^\flat\nu^\times) / \kappa^+               &
                                                               \sigma_{21}^\times + \sigma_2^+ - \nu^+   \\
  -(1+\eps) \sigma_{12}^\times \kappa^\times / \kappa^+  + \sigma_1^+ - \kappa^+ &        0              \\
  \sigma_{22}^\times (\chi\nu^+ + \eps^\flat \nu^\times) / \kappa^+              &
                                                               \sigma_{22}^\times + \sigma_2^+ - \nu^+
           \end{pmatrix}.
\end{eqnarray*}

\begin{prop}
\label{prop:ker}
The vector $\begin{pmatrix} a \\ b \end{pmatrix}$,
$a = \begin{pmatrix} 0 \\ \kappa^- \end{pmatrix}$,
$b = \begin{pmatrix} \kappa^+ \\ (\eps^\times - \chi) \nu^+ \end{pmatrix}$
lies in the kernel of $\Theta|_{\lambda=0}$.
\end{prop}

This can be checked by a direct computation. However this is a general fact, see Remark \ref{rem:lam=0}.

\subsubsection{Completion of the proof of Theorem \ref{theo:main}}

By (\ref{sigma-}), (\ref{sigma+}), and (\ref{sigmatimes}) we obtain $\Theta = \Theta_1 + \Theta_2$, where
\begin{eqnarray*}
    \Theta_1
&=& \begin{pmatrix}
     -\nu^-  &  \frac{\eps\lambda}{\kappa^-}       & -\frac{2+\eps}{\kappa^+} \lambda  &       0  \\[1mm]
       0     &  \frac{(2+\eps)\lambda}{\kappa^-}   &
              \frac{\chi\nu^+ + \eps^\flat\nu^\times}{\nu^\times \kappa^+} \lambda     &
                                        \big(\frac1{\nu^\times} - \frac1{\nu^+}\big) \lambda      \\[1mm]
   \frac{1+\eps}{\eps} (\kappa^\times - \nu^-)
             &  -\frac{\kappa^{\times\, 2}(1+\eps)}{\kappa^-}
                                                   & \frac{(1+\eps)\kappa^{\times\, 2}}{\kappa^+} &  0 \\[1mm]
        0    &    -\frac{\nu^{\times\, 2}(\eps - \eps^\times)}{\kappa^-} + \frac\lambda{\kappa^-}
                         &  - \frac{(\chi\nu^+ + \eps^\flat \nu^\times)}{\kappa^+} \nu^\times
                                                                 & - \nu^\times
     \end{pmatrix}, \\
    \Theta_2
&=& \begin{pmatrix}
     O(\lambda |\log\eps|)     &  O(\eps\lambda^2|\log\eps|)         &   O(\lambda^2) &        0       \\
       0                       &  O(\lambda^2|\log\eps|)             &   O(\lambda^2) &  O(\lambda^2)  \\
 O(\lambda\eps^{-1}|\log\eps|) &  O(\lambda\log\eps)                 &   O(\lambda)   &        0       \\
       0                       & (\eps-\eps^\times) O(\lambda|\log\eps|) + O(\lambda^2)
                                                                     &  O(\lambda)   & O(\lambda)
     \end{pmatrix}.
\end{eqnarray*}

\begin{prop}
\label{prop:d}
Suppose $2\nu^\times + \chi < 0$. Then for sufficiently small $\eps > 0$ the equation   $\det\Theta=0$ has a positive solution
\begin{equation}
\label{lambda(eps)}
    \lambda(\eps)
  = \frac{\nu^- {\kappa^\times}^2 (2\nu^\times + \chi)}{4\nu^\times(\kappa^\times - \nu^-)}\eps + o(\eps).
\end{equation}
\end{prop}

\noindent
{\it Proof}. Let
$$
    \bd(\eps,\lambda)
  = \frac{\kappa^+ \kappa^-}{(1+\eps)\nu^\times\lambda} \det\Theta.
$$
Then $\bd = \bd_1 + \bd_2$, where
$$
     \bd_1
  =  \Big( -\nu^- + O(\lambda |\log\eps|) \Big) \det(\Phi_1 + \Psi_1), \quad
     \bd_2
  =  \Big( \frac{\kappa^\times - \nu^-}{\eps} + O(\lambda\eps^{-1} |\log\eps|) \Big)
       \det(\Phi_2 + \Psi_2).
$$
The matrices $\Phi_1,\Phi_2,\Psi_1$,$\Psi_2$ are as follows:
\begin{eqnarray*}
     \Phi_1
 &=& \begin{pmatrix}
   2+\eps           & \frac{\chi\nu^+ + \eps^\flat\nu^\times}{\nu^\times} & \frac1{\nu^\times}-\frac1{\nu^+} \\
 -\kappa^{\times\, 2} & \kappa^{\times\, 2}                               &     0 \\
 -\nu^\times (\eps-\eps^\times) + \frac\lambda{\nu^\times}
                      & -(\chi\nu^+ + \eps^\flat\nu^\times)               &  -1
        \end{pmatrix} , \\
     \Psi_1
 &=& \begin{pmatrix}
       O(\lambda|\log\eps|)                                 &   O(\lambda)   &  O(\lambda)  \\
       O(\lambda\log\eps)                                   &   O(\lambda)   &        0       \\
     (\eps-\eps^\times) O(\lambda|\log\eps|) + O(\lambda^2) &   O(\lambda)   & O(\lambda)
     \end{pmatrix}, \\
     \Phi_2
 &=&
 \begin{pmatrix}
   \eps\lambda      & -(2+\eps)\lambda                                    &     0  \\
   2+\eps           & \frac{\chi\nu^+ + \eps^\flat\nu^\times}{\nu^\times} & \frac1{\nu^\times}-\frac1{\nu^+} \\
 -\nu^\times (\eps-\eps^\times) + \frac\lambda{\nu^\times}
                      & -(\chi\nu^+ + \eps^\flat\nu^\times)               &  -1
        \end{pmatrix},
        \\
     \Psi_2
 &=&
 \begin{pmatrix}
             O(\eps\lambda^2|\log\eps|)                  &   O(\lambda^2) &        0     \\
             O(\lambda|\log\eps|)                        &   O(\lambda)   &  O(\lambda)  \\
 (\eps-\eps^\times) O(\lambda|\log\eps|) + O(\lambda^2)  &  O(\lambda)    &  O(\lambda)
\end{pmatrix}.
\end{eqnarray*}
We obtain by a direct computation:
$$
      \bd_1
   =  \nu^-\kappa^{\times\, 2} \Big(2 + \frac\chi{\nu^\times}\Big) + \calO, \quad\!
      \bd_2
   =  -4\lambda\frac{\kappa^\times - \nu^-}{\eps} (1 + \calO), \qquad\!\!
      \calO = O(\eps + \eps^\times + \eps^\flat) + O(\lambda |\log\eps|).
$$
Therefore by the implicit function theorem the equation $\bd=0$ has a solution of the form (\ref{lambda(eps)}).
\qed

\medskip

Now Theorem \ref{theo:main} is proved.

\section{Appendix}

\label{sec:App}

\subsection{Quadratic inequalities}

\label{sec:quad}

For any real $h\ne 0$ we define
\begin{equation}
\label{P_h}
  \calP_h^- = \{\lambda\in\mC : h^2 \re\lambda < - (\im\lambda)^2\},\quad
  \calP_h^+ = \{\lambda\in\mC : h^2 \re\lambda > - (\im\lambda)^2\}.
\end{equation}

\begin{lem}
\label{lem:h<0}
{\bf (1)}. Suppose $h<0$. Then the polynomial $p_h(\sigma) = \sigma^2 - h\sigma - \lambda$ has
\begin{itemize}
\item no roots with $\re\sigma > 0$ for $\lambda\in \calP^-_h$,
\item one root with $\re\sigma > 0$ for $\lambda\in \calP^+_h$.
\end{itemize}

{\bf (2)}. Suppose $h>0$. Then the polynomial $p_h(\sigma) = \sigma^2 - h\sigma - \lambda$ has
\begin{itemize}
\item one root with $\re\sigma > 0$ for $\lambda\in \calP^+_h$,
\item two roots with $\re\sigma > 0$ for $\lambda\in \calP^-_h$.
\end{itemize}
\end{lem}

\noindent
{\it Proof}. {\bf (1)}. Sum of the roots equals $h < 0$. Therefore at most one root has positive real part. Suppose such a root $\sigma = \alpha + i\beta$ exists. Then
$$
  \lambda = (\alpha + i\beta)^2 - h(\alpha + i\beta)
          = (\alpha + i\beta - h/2)^2 - h^2/4.
$$
The domain
$\{\lambda\in\mC : \lambda = (\alpha + i\beta - h/2)^2 - h^2/4, \quad \alpha > 0,\; \beta\in\mR \}$
coincides with $\calP_h^+$.

{\bf (2)}. Sum of the roots equals $h > 0$. Therefore at least one root has positive real part. Suppose there are two such roots. This happens if both of them have the form $\sigma = \alpha + i\beta$, $0 < \alpha < h$. Then $\lambda$ belongs to the domain
$$
  \{\lambda\in\mC : \lambda = (\alpha + i\beta - h/2)^2 - h^2/4, \quad 0 < \alpha < h, \; \beta\in\mR \}.
$$
This domain coincides with $\calP_h^-$. \qed

\begin{lem}
\label{lem:h>0}
{\bf (1)}. Suppose $h<0$. Then the polynomial $p_h(\sigma) = \sigma^2 - h\sigma - \lambda$ has
\begin{itemize}
\item one root with $\re\sigma < 0$ for $\lambda\in \calP^+_h$,
\item two roots with $\re\sigma < 0$ for $\lambda\in \calP^-_h$.
\end{itemize}

{\bf (2)}. Suppose $h>0$. Then the polynomial $p_h(\sigma) = \sigma^2 - h\sigma - \lambda$ has
\begin{itemize}
\item no roots with $\re\sigma < 0$ for $\lambda\in \calP^-_h$,
\item one root with $\re\sigma > 0$ for $\lambda\in \calP^+_h$.
\end{itemize}
\end{lem}

\noindent
{\it Proof}. If $\sigma_1,\sigma_2$ are roots of $p_h$ then $p_{-h}$ has the roots $-\sigma_1,-\sigma_2$. This reduces Lemma \ref{lem:h>0} to Lemma \ref{lem:h<0}. \qed

\subsection{Proof of Proposition \ref{prop:smooth}}

\label{sec:smoothing}

For simplicity we write a vector $u\in\mR^n$ as
$u=(u_1,u_2)$, $u_1\in\mR^{n-1}$, $u_2\in \mR$, and let  $\Sigma=\{u_2=0\}$.
Then $U_\mu=\{0<u_2<\mu\}$ and equation (\ref{eq:phi}) becomes
$$
f^\mu(u)=\phi(u_2/\mu)f^-(u)+(1-\phi(u_2/\mu))f^+(u).
$$
First consider the trajectory $\gamma^\mu(\xi)$ of $\gamma'=f^\mu(\gamma)$
crossing $\Sigma$ at $\gamma^\mu(0)=w=(w_1,0)$.

\begin{lem}
$\gamma^\mu(\xi)$ crosses the layer $U_\mu$ for $0<\xi<\mu a$, where
$$
a(\mu)=\int_0^1\frac{dy}{\phi(y)\hat f_2^++(1-\phi(y))\hat f_2^-}+O(\mu),\qquad \hat f^\pm=f(w_1,0).
$$
For $0<\xi<\mu a$, we have
$$
\gamma^\mu(\xi)=(w_1^0+O(\mu),\mu\Phi^{-1}(\mu\xi)+O(\mu^2)),\quad
\Phi(y)=\int_0^y\frac{dx}{\phi(x)\hat f_2^++(1-\phi(x))\hat f_2^-}
$$
\end{lem}

\noindent{\it Proof}. After the change of variables
$$
u_2=\mu y,\qquad 0<y<1,\qquad  \xi=\mu s,
$$
in $U_\mu$, the system $u'=f^\mu(u)$ becomes
\begin{eqnarray*}
\frac{ d u_1}{ds}&=&\mu \phi(y)f_1^+(u_1,\mu y)+\mu (1-\phi(y))f_1^-(u_1,\mu y)=O(\mu),\\
\frac{dy}{ds}&=& \phi(y)f_2^+(u_1,\mu y)+(1-\phi(y))f_2^-(u_1,\mu y)=\phi(y)f_2^+(u_1,0)+(1-\phi(y))f_2^-(u_1,0)+O(\mu).
\end{eqnarray*}
Hence on a finite $s$-interval, $u_1=w_1+O(\mu)$ and
$$
\frac{dy}{ds}=\phi(y)\hat f_2^++(1-\phi(y))\hat f_2^-+O(\mu),\qquad \hat f_i^\pm=f_i^\pm(w_1,0),\quad y(0)=0.
$$
It remains to separate the variables.
\qed

\medskip

 \noindent
 {\it Proof of Proposition \ref{prop:smooth}}. Consider the smoothed   spectral problem (\ref{eq:zv}).
 We have $\partial f^\mu=(\partial_1 f^\mu,\partial_2 f^\mu)$, where
 \begin{eqnarray*}
\partial_1 f^\mu&=&
 \phi(y)\partial_1 f^+(u_1,\mu y)+(1-\phi(y))\partial_1 f^-(u_1,\mu y), \\
\partial_2 f^\mu&=&\phi(y)\partial_2 f^+(u_1,\mu y)+(1-\phi(y))\partial_2 f^-(u_1,\mu y)
 +\mu^{-1}\phi'(y)(  f^+(u_1,\mu y) -  f^-(u_1,\mu y))
\end{eqnarray*}
Let $v=(v_1,v_2)$. System (\ref{eq:zv}) takes the form
\begin{eqnarray*}
\frac{ dz}{ds}&=&\mu\lambda v,\qquad 0\le s\le a,\\
\frac{ dv}{ds}&=&\mu z+\mu\big(\phi(y)\partial_1 f^+(u_1,\mu y)+(1-\phi(y))\partial_1 f^-(u_1,\mu y)\big)v_1\\
&&\qquad +\mu \big(\phi(y)\partial_2 f^+(u_1,\mu y)+(1-\phi(y))\partial_2 f^-(u_1,\mu y)\big) v_2\\
&&\qquad +\phi'(y)(  f^+(u_1,\mu y) -  f^-(u_1,\mu y))v_2\\
&=&\phi'(y)(  \hat f^+ -  \hat f^-)v_2+O(\mu).
\end{eqnarray*}

Hence  for $0\le s\le a$,
$$
z(s)=z(0)+O(\mu),\quad   \frac{dv}{ds}=\phi'(y)\Delta f v_2+O(\mu).
$$
We use  $y$ as an independent variable:
$$
  \frac{dv}{dy} = \frac{\phi'(y)\Delta f v_2}{\phi(y)\hat f_2^+ +(1-\phi(y))\hat f_2^-}+O(\mu),\qquad
  0\le y\le 1,
$$
Then pass to the independent variable $\phi$:
$$
  \frac{dv}{d\phi} = \frac{\Delta f  v_2}{\phi \hat f_2^+ +(1-\phi)\hat f_2^-} +O(\mu),\qquad
  0\le \phi\le 1.
$$
In the first approximation
\begin{equation}
\label{eq:v12}
  \frac{dv_1}{d\phi} = \frac{\Delta f_1\, v_2}{\phi \hat f_2^+ +(1-\phi)\hat f_2^-},\qquad
  \frac{dv_2}{d\phi} = \frac{\Delta f_2\, v_2}{\phi \hat f_2^+ +(1-\phi)\hat f_2^-}  .
\end{equation}
First solve the second equation:
$$
    \ln v_2-\ln v_2(0)=\ln (\phi \hat f_2^+ +(1-\phi)\hat f_2^-)-\ln f_2^-,\quad
            v_2=\frac{\phi \hat f_2^+ +(1-\phi)\hat f_2^-}{\hat f_2^-} v_2(0).
$$
 From the first equation (\ref{eq:v12}),
$$
  v_1-v_1(0) = \phi \frac{\Delta f_1} {\hat f_2^-}v_2(0) ,\qquad 0\le \phi\le 1.
$$
Let $v^-=v(0)$, $v^+=v|_{\phi=1}$, then
$$
\label{eq:Delta_v}
   v_2^+ = \frac{\hat f_2^+ }{\hat f_2^-} v_2^-,\quad
   v_1^+ = v_1^- + \frac{\Delta f_1 }{\hat f_2^-}v_2^-,\quad
   v^+ = v^- + \frac{\Delta f }{\hat f_2^-}v_2^-.
$$
Proposition \ref{prop:smooth} is proved.
\qed

\section*{Acknowledges}

The work of S.~Bolotin  work was performed at the Steklov International Mathematical Center
and supported by the Ministry of Science and Higher Education of the Russian Federation
(agreement no.\ 075-15-2025-303).
The work of D.~Treschev was supported by the Russian Science Foundation under grant no.\ 25-11-00114.
S.~Bolotin worked on Sections 2, 4, 6, and D.~Treschev on Sections 1, 3,  5.


\begin{thebibliography}{00}

\bibitem{Chicone}
C. Chicone,
Quadratic gradients on the plane are generically Morse-Smale. J. Differ. Eq.
 {\bf 33} (1979), 159--166.

\bibitem{Chug-Tre}
A. Chugainova and D. Treschev,
Waves, structures, and the Riemann problem
for a hyperbolic system of conservation laws.
arXive http://arxiv.org/abs/2510.02070.

\bibitem{Chug-Pol}
A.P. Chugainova and R.R. Polekhina, Admissibility of discontinuities in the solutions of a
hyperbolic $2\times 2$ system of conservation laws,
International Journal of Non-Linear Mechanics,  {\bf 178} (2025), 105174.

\bibitem{Ilich}
V.I. Erofeev and A.T. Il'ichev,
Instability of supersonic solitary waves in a generalized elastic electrically conductive medium,
Contin. Mech. Thermodin., {\bf 35} (2023), 2313--2323.

\bibitem{Evans}
J.W. Evans, Nerve axon equations: I. Linear approximations. Ind. Univ.
Math. J. {\bf 21} (1972), 877--885.

\bibitem{Gar-Zum}
R. A. Gardner and K. Zumbrun,
The gap lemma and geometric criteria for instability of viscous shock profiles,
Comm. Pure Appl. Math., {\bf 51} (1998), 797--855.

\bibitem{Good}
J. Goodman,
Nonlinear asymptotic stability of viscous shock profiles for conservation laws,
Arch. Ration. Mech. Anal., {\bf 95} (1986), 325--344.

\bibitem{How-Zum}
P. Howard and K. Zumbrun,
Stability of undercompressive shock profiles,
J. Differential Equations {\bf 225} (2006), 308--360.

\bibitem{IO}
A.M. Il'in and O.A. Oleinik,
Asymptotic behavior of solutions of the Cauchy problem for some quasi-linear equations
for large values of the time. Mat. Sb. (N.S.) 51(93) (1960), 191--216.

\bibitem{Liu}
T.-P. Liu,
Pointwise convergence to shock waves for viscous conservation laws.
Comm. Pure Appl. Math. {\bf 50} (1997), 1113--1182.

\bibitem{LZ1}
T.-P. Liu and K. Zumbrun,
Nonlinear stability of an undercompressive shock for complex Burgers equation.
Comm. Math. Phys., {\bf 168} (1995), 163--186.

\bibitem{LZ2}
T.-P. Liu and K. Zumbrun,
On nonlinear stability of general undercompressive viscous shock waves.
Comm. Math. Phys., {\bf 174} (1995), 319--345.

\bibitem{MN}
A. Matsumura and K. Nishihara,
On the stability of travelling wave solutions of a one-dimensional model system for compressible viscous gas.
Japan J. Appl. Math., {\bf 2} (1985), 17--25.

\bibitem{Satt}
D.H. Sattinger,
On the stability of waves of nonlinear parabolic systems.
Adv. Math., {\bf 22} (1976), 312--355.

\bibitem{Zum-How}
K. Zumbrun and P. Howard,
Pointwise semigroup methods and stability of viscous shock waves,
Indiana Univ. Math. J., {\bf 47} (1998), 741--871.

\bibitem{Zum:Survey}
K. Zumbrun,
Stability and dynamics of viscous shock waves.
In: Bressan, A., Chen, GQ., Lewicka, M., Wang, D. (eds)
Nonlinear Conservation Laws and Applications.
The IMA Volumes in Mathematics and its Applications, vol 153.
Springer, Boston, MA. (2011)

\end{thebibliography}
\end{document}